\documentclass[final]{siamltex}

\usepackage{amsmath}
\usepackage{graphicx}
\usepackage{setspace}
\usepackage{citesort}
\usepackage{SIunits}
\usepackage{url}

\newcommand{\bs}{\hat{\mathbf{s}}}
\newcommand{\bw}{\boldsymbol{\omega}}

\newcommand{\D}{\mathrm{d}}

\begin{document}

\title{Evaluation of Biot-Savart integrals on tetrahedral
  meshes\thanks{This work was supported by the European
    Community-funded project HPC-Europa, contract number 506079, and
    was carried out at the Institut f\"{u}r Str\"{o}mungsmechanik und
    Hydraulische Str\"{o}mungsmaschine and HLRS, both at the
    University of Stuttgart, Germany}}

\author{Michael Carley\thanks{Department of Mechanical
    Engineering, University of Bath, Bath BA2 7AY, United Kingdom
    ({\tt{m.j.carley@bath.ac.uk})}}}

\bibliographystyle{siam}

\maketitle

\begin{abstract}
  An arithmetically simple method has been developed for the
  evaluation of Biot--Savart integrals on tetrahedralized
  distributions of vorticity. In place of the usual approach of
  analytical formulae for the velocity induced by a linear
  distribution of vorticity on a tetrahedron, the integration is
  performed using Gaussian quadrature and a ray tracing technique from
  computer graphics. This eliminates completely the need for the
  evaluation of square roots, logarithms and arc tangents, and almost
  completely eliminates the requirement for trigonometric functions,
  with no operation more costly than a division required during the
  main calculation loop. An assessment of the algorithm's performance
  is presented, demonstrating its accuracy, second order convergence
  and near-linear speedup on parallel systems.
\end{abstract}



\pagestyle{myheadings}

\thispagestyle{plain}

\markboth{M. CARLEY}{Biot-Savart integrals}

\section{Introduction}
\label{sec:intro}

An important part of many calculations in fluid dynamics and
electromagnetism is the evaluation of a Biot-Savart integral, for
velocity in fluid dynamics and for magnetic field in electromagnetism.
In fluid dynamics, the source term is vorticity while in
electromagnetism it is current. The Biot-Savart integral for velocity
$\mathbf{v}$ due to a distribution of vorticity $\bw$ over a volume
$V$ is:
\begin{align}
  \label{equ:biot:savart:basic}
  \mathbf{v}(\mathbf{x}) &= 
  -\frac{1}{4\pi}\int_{V}
  \frac{\mathbf{r}\times\bw(\mathbf{x}_{1})}{R^{3}}\,\D V,
\end{align}
where $\mathbf{r}=\mathbf{x}-\mathbf{x}_{1}$, $R=|\mathbf{r}|$ and
subscript~1 indicates a variable of integration. For electromagnetic
calculations, $\mathbf{v}$ is the magnetic field and $\bw$ the current
in the region $V$. For the purposes of this paper, it is assumed that
the volume $V$ is discretized into tetrahedral elements within which
the source varies linearly. The question then is how to compute the
resulting field $\mathbf{v}$. This paper is motivated by the
requirement to evaluate velocities in Lagrangian vortex methods, where
a moving distribution of control points is tetrahedralized at each
time step as an aid to velocity calculation. Such a method has been
implemented by Marshall et al.~\cite{marshall-grant-gossler-huyer00},
using a mixture of analytical formulae and Gaussian quadrature to
carry out the integration of equation~\ref{equ:biot:savart:basic}. The
aim of this paper is to develop a method which is arithmetically
simpler than that used hereto. 

A number of analytical formulae and numerical procedures have been
developed for the evaluation of equation~\ref{equ:biot:savart:basic}.
A recent general paper is that of Suh~\cite{suh00} who applies Stokes'
theorem to reduce the volume integral over an element to a number of
surface integrals which are in turn reduced to line integrals which
can be evaluated analytically. This work is similar to that of
Newman~\cite{newman86} who derived equivalent formulae, focussing on
applications in fluid dynamics. In electromagnetism, there is an
extensive literature on the evaluation of
equation~\ref{equ:biot:savart:basic}, including
analytical~\cite{urankar84iv,weggel-schwartz88,onuki-wakao95,graglia87,%
  wilton-rao-glisson-schaubert-al-bundak-butler84} and
numerical~\cite{khayat-wilton05} approaches. As this paper is
motivated by the fluid dynamical problem, it will use the terms
`velocity' and `vorticity' but it should be noted that there are also
many useful related results in the electromagnetism literature.

In evaluating the velocity at the nodes of a vorticity distribution
which is discretized into tetrahedra, a Gaussian quadrature can be
used for tetrahedra which are far from the evaluation point
$\mathbf{x}$, but (integrable) singularities in the integrand make
this awkward for nearby elements. In this case, the standard approach
is to use analytical formulae which have the benefit of being exact
and non-singular. The disadvantage, as noted by Marshall et
al.~\cite{marshall-grant-gossler-huyer00}, is that such
formulae~\cite{newman86} are computationally expensive, requiring in
this case~12 logarithms and~24 arc tangents per tetrahedron. Even if
the analytical formulae are only used for `nearby' elements, they
still represent a large part of the velocity computation. The aim of
this paper is to present a velocity computation method which is
arithmetically simple, requiring no operation more onerous than a
division during the main calculation, by using a ray-tracing method
borrowed from computer graphics. The resulting method can then be used
in codes as a `plug-in' replacement for previous techniques.

\section{Velocity evaluation}
\label{sec:velocity}

The integral of equation~\ref{equ:biot:savart:basic} is to be
evaluated at a number of points $\mathbf{x}$, which also form the
nodes of a distribution of vorticity $\bw$. In this case, as in the
work of Marshall et al.~\cite{marshall-grant-gossler-huyer00}, the
nodes are tetrahedralized so that they form the vertices of a
collection of tetrahedra. It is further assumed that vorticity varies
linearly over the elements. The aim is now to avoid the computational
effort involved in the standard analytical approach to quadrature and
develop a method which is arithmetically as simple as possible. In the
words of Richardson, this is a problem where it may be quicker to
arrive at a destination at the ``footpace of arithmetic'', rather than
``on the swift steed of symbolic analysis''~\cite{richardson25}.

First, equation~\ref{equ:biot:savart:basic} is rewritten in spherical
polar coordinates centred on the evaluation point $\mathbf{x}$:
\begin{align}
  \label{equ:biot:savart}
  \mathbf{v} &= \frac{1}{4\pi}
  \int_{0}^{\pi}
  \int_{0}^{\pi}
  \int_{-\infty}^{\infty}
  \bs\times\bw
  \,\D R
  \sin\phi\,\D\phi
  \D\theta,
\end{align}
where $\mathbf{r}=\mathbf{x}-\mathbf{x}_{1}=-R\bs$ with the unit
vector $\bs=(\sin\phi\cos\theta,\sin\phi\sin\theta,\cos\phi)$. This
transformation also eliminates the integrable singularity in the
integrand, making the integration rather easier from a numerical point
of view. The azimuthal and polar integrals are evaluated using
Gaussian quadratures so that:
\begin{align}
  \label{equ:biot:savart:num}
  \mathbf{v} &\approx 
  \sum_{n=1}^{N}\sum_{m=1}^{M}
  \sin\phi_{n}\frac{w^{(N)}_{n}w^{(M)}_{m}}{16\pi^{2}}
  \int_{-\infty}^{\infty}
  \bs_{nm}\times\bw
  \,\D R,
\end{align}
where
\begin{align*}
  \phi_{n} &= (1+t^{(N)}_{n})\pi/2,\\
  \theta_{m} &= (1+t^{(M)}_{m})\pi/2,\\
  \bs_{nm} &=
  (\sin\phi_{n}\cos\theta_{m},\sin\phi_{n}\sin\theta_{m},\cos\phi_{n}),
\end{align*}
and $(t_{i}^{(N)},w_{i}^{(N)})$, $i=1,\ldots,N$, are the nodes and
weights of an $N$ point Gaussian quadrature with $-1\leq t\leq 1$. The
main computation is now reduced to integrals over $R$ along rays in
the direction $\bs$. 


\subsection{Integration over tetrahedra}
\label{sec:tetrahedra}

\begin{figure}
  \centering
  \includegraphics{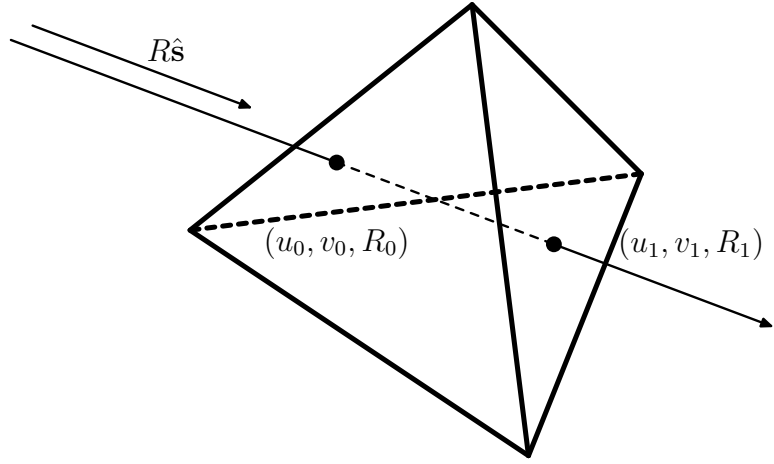}
  \caption{Intersection of a ray with a tetrahedron: the ray from the
    evaluation point in direction $\bs$ intersects two faces of the
    tetrahedron at the points marked by circles}
  \label{fig:intersection}
\end{figure}

The radial integral of equation~\ref{equ:biot:savart:num} must now be
evaluated, to include the contribution of the mesh tetrahedra. The
approach is summarized in figure~\ref{fig:intersection}: if the ray
from $\mathbf{x}$ in the direction $\bs$ intersects the tetrahedron,
i.e. passes through two of its faces, the tetrahedron makes a finite
contribution to the integral. The overall integral is:
\begin{align}
  \label{equ:ray:int}
  I(\mathbf{x},\theta,\phi) &= 
  \int_{-\infty}^{\infty}
  \bs\times\bw
  \,\D R.
\end{align}
If the ray enters a tetrahedron at $R=R_{0}$ and exits at $R=R_{1}$,
\begin{align}
  \label{equ:ray:sum}
  I(\mathbf{x},\theta,\phi) &= 
  \sum   \int_{R_{0}}^{R_{1}}
  \bs\times\bw
  \,\D R &=
  \frac{1}{2}\sum
  (R_{1}-R_{0})\bs\times(\bw_{0}+\bw_{1}), 
\end{align}
summing over tetrahedra which are intersected by the ray. The
vorticities $\bw_{0}$ and $\bw_{1}$ are those at the entry and exit
points $R_{0}$ and $R_{1}$ and, as noted previously, $\bw$ varies
linearly between them, so that equation~\ref{equ:ray:sum} is exact. 
If the ray enters or exits a face of the tetrahedron at area
coordinates $(u,v)$, the vorticity $\bw$ is:
\begin{align*}
  \bw &= (1-u-v)\bw_{0} + u\bw_{1} + v\bw_{2},
\end{align*}
where $\bw_{i}$, $i=0,1,2$ are the vorticities at the triangle
vertices. This integral is arithmetically very simple, with no
operation more complex than a division, not even requiring a square
root for the calculation of a distance. 

The remaining part of the algorithm is how the intersection, if any,
of a ray with a tetrahedron can be determined. This is found using a
method from computer graphics. M\"{o}ller and
Trumbore~\cite{moller-trumbore97} give a method and C source code for
determining the area coordinates $(u,v)$ and directed distance $R$ of
the intersection point of a ray with a triangle. Their method is very
efficient requiring at most one division and no operation more complex
than this. It is used in the algorithm of this paper to find the
intersections of a ray with a tetrahedron, by checking the faces in
turn. In order to accelerate the procedure, an initial check is
carried out on the tetrahedron as a whole.

\begin{figure}
  \centering
  \includegraphics{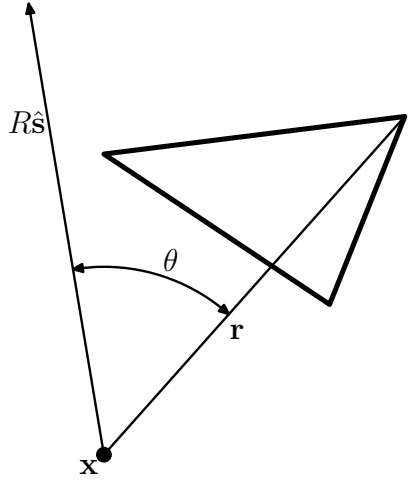}
  \caption{Pre-check to decide if tetrahedron is likely to be
    intersected by a ray} 
  \label{fig:precheck}
\end{figure}

Figure~\ref{fig:precheck} shows the check for the case of a
triangle. The vector $\mathbf{r}$ is the vector from $\mathbf{x}$ to
an arbitrary vertex of the tetrahedron. The scalar product of $\bs$
and $\mathbf{r}$ is:
\begin{align*}
  \mathbf{r}.\bs = |\mathbf{r}|\cos\theta,
\end{align*}
where $\theta$ is the angle between $\bs$ and $\mathbf{r}$. If the
tetrahedron subtends less than this angle, none of its faces are
intersected and no further check is needed. If it has a maximum edge
length $h$, the maximum angle it can subtend is $\psi\approx
h/|\mathbf{r}|$ and $\cos\theta<\cos\psi$ is a necessary condition for
the tetrahedron to be intersected. Using the small angle approximation
of $\cos\psi$ and neglecting fourth order terms:
\begin{align*}
  \cos\theta < 1-h^{2}/|\mathbf{r}|^{2},\\
  \cos^{2}\theta < 1-2h^{4}/|\mathbf{r}|^{4},
\end{align*}
so that the condition can be written:
\begin{align*}
  (\mathbf{r}.\bs)^{2} < |\mathbf{r}|^{2} - 2h^{2}.
\end{align*}
The test is implemented by setting $h^{2}$ to the square of the
longest edge length on the element and by applying it only for
$|\mathbf{r}|$ greater than some minimum value.

\subsection{Summary of algorithm}
\label{sec:summary}

In summary, given a set of points $\mathbf{x}_{i}$ each with an
associated vorticity $\bw_{i}$, the Biot-Savart integral at each point
can be evaluated as follows:
\begin{enumerate}
\item Generate a tetrahedralization of the points $\mathbf{x}_{i}$
  using, for example, the TETGEN code of Si~\cite{si06}.
\item Select Gaussian quadrature rules of order $N$ and $M$ for
  integration in $\phi$ and $\theta$ respectively.
\item Precompute the factors
  $k_{nm}=\sin\phi_{n}w_{n}^{(N)}w_{m}^{(M)}/16\pi^{2}$ and the ray
  directions $\bs_{nm} =
  (\sin\phi_{n}\cos\theta_{m},\sin\phi_{n}\sin\theta_{m},\cos\phi_{n})$,
  equation~\ref{equ:biot:savart:num}, for $n=1,\ldots,N$ and
  $m=1,\ldots,M$.
\item For each tetrahedron, loop over the nodes $\mathbf{x}_{i}$ and
  check for intersection of each ray $\bs$ in turn. If the ray
  intersects the tetrahedron, add the contribution of
  equation~\ref{equ:ray:sum} to the velocity $\mathbf{v}_{i}$.
\end{enumerate}
If the calculation is being carried out on a parallel system, the
tetrahedra can be divided amongst the processors for the velocity
calculation with the total velocity being found by a summation at the
end of the computation. This only requires one communication per
velocity calculation. 

\section{Numerical tests}
\label{sec:tests}

A number of tests have been carried out to evaluate the performance of
the method of \S\ref{sec:summary} with respect to speed, accuracy and
parallelization. The method has been implemented in a Lagrangian
vortex code currently under development, similar to that of Marshall
et al.~\cite{marshall-grant-gossler-huyer00}. The control points are
tetrahedralized using the TETGEN code of Si~\cite{si06} and the
program runs on serial and parallel systems using the MPI message
passing interface. 

\subsection{Speed and accuracy}
\label{sec:speed}

\begin{table}
  \caption{Computation time and r.m.s. error in velocity calculations
    for Hill's spherical vortex}
  \label{tab:performance}
  \centering
  \begin{tabular}{rrrrrrr}
    &
    \multicolumn{2}{c}{Summation~\cite{marshall-grant-gossler-huyer00}} &
    \multicolumn{4}{c}{Ray tracing} \\
    $N$
    &  
    \multicolumn{1}{c}{$T/\second$}
    & 
    \multicolumn{1}{c}{$\epsilon$}
    &  
    \multicolumn{1}{c}{$T/\second$}
    & 
    \multicolumn{1}{c}{$\epsilon$}
    &  
    \multicolumn{1}{c}{$T/\second$}
    & 
    \multicolumn{1}{c}{$\epsilon$}\\
    \hline
    1000 & 2 & 0.200 & 28 & 0.036 & 13 & 0.036\\
    2000 & 8 & 0.100 & 100 & 0.028 & 44 & 0.026\\
    4000 & 30 & 0.020 & 369 & 0.021 & 155 & 0.017\\
    8000 & 100 & 0.006 & 1385 & 0.019 & 562 & 0.013\\
    16000 & 400 & 0.004 &  &  & 2168 & 0.009\\
    \hline
    & \multicolumn{2}{c}{Cray C90}
    & \multicolumn{2}{c}{Laptop}
    & \multicolumn{2}{c}{8-node cluster}
  \end{tabular}
\end{table}

The first test case considered is Hill's spherical
vortex~\cite[pp23--25]{saffman92}. This is a distribution of vorticity
with a linear variation of azimuthal vorticity $\bw_{\theta}=Ar$
inside a sphere of radius $a$. Within the sphere, radial and axial
velocities are:
\begin{align}
  \label{equ:hill}
  u &= \frac{Az}{5},\quad
  v = \frac{A}{5}\left(2r^{2} + z^{2} + \frac{5}{3}a^{2}\right). 
\end{align}
This is an especially useful test case because, with the exception of
errors in discretizing the boundary of the sphere, the vorticity
distribution is exactly duplicated by linear interpolation over
elements. The test was carried out by evaluating the velocity at $N$
randomly placed points with a sphere of radius $a=1$ with $N$ varying
from 1000 to 16000. Calculations were carried out a 500\mega\hertz\
Intel Pentium III laptop with~4 point Gaussian quadrature in $\phi$
and $\theta$ and with~16 point quadratures on an~8 node
3.2\giga\hertz\ Intel Xeon cluster. Table~\ref{tab:performance} shows
the r.m.s error in velocity and the computation time for these two
cases and, for comparison, the corresponding data for Marshall et
al.'s calculations on a Cray
C-90~\cite{marshall-grant-gossler-huyer00}. From the data presented,
it appears that for small numbers of nodes, the method is more
accurate than the use of analytical formulae and that it remains
better up to distributions, in this case, of about~4000 points. The
second point to note is that the accuracy appears to be controlled, up
to the limit of the resolution of the vorticity, by the number of
quadrature points.  The $4\times4$ quadrature error stops reducing at
about~4000 nodes while the $16\times16$ is still improving at~16000
nodes.

\subsection{Convergence}
\label{sec:convergence}

\begin{figure}
  \centering
  \includegraphics{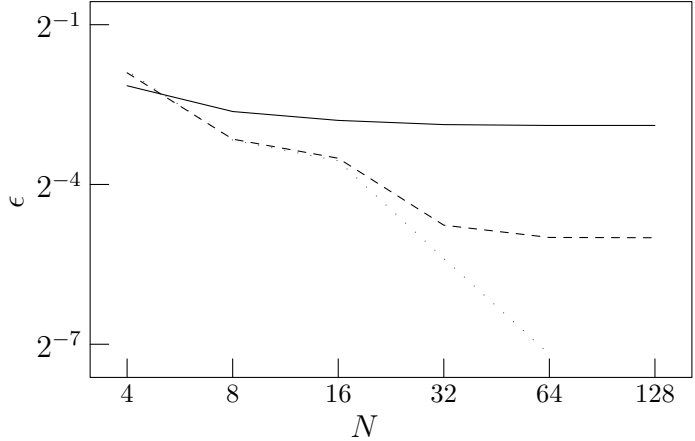}
  \caption{Convergence of velocity for a Gaussian vortex ring, r.m.s.
    error against number of quadrature points: solid line $16\times17$
    points; dashed line $64\times65$ points; dotted line
    $128\times257$ points.}
  \label{fig:convergence}
\end{figure}

As a second test on the accuracy of the calculation method, and to
assess the convergence properties, the velocity due to a Gaussian core
ring was computed. This has exponentially small vorticity on the
boundary of the vorticity distribution so that errors due to
discretization of the boundary should be reduced, compared to the case
of Hill's spherical vortex. For comparison, the velocity was computed
using the stream function for an axisymmetric vortex
ring~\cite[pages~192--194]{saffman92}, differentiating it numerically
on a dense grid to give the velocity. This velocity was interpolated
to find the velocity at the control points of the test
distribution. Three different distributions of points were considered,
each made up of a number of stations equally spaced in azimuth, with a
square grid of points at each station. The low resolution distribution
had~$16\times17$ points, the intermediate~$64\times65$ and the high
resolution~$128\times257$. The velocity was evaluated using equal
numbers of quadrature points in $\phi$ and
$\theta$. Figure~\ref{fig:convergence} shows the r.m.s. error in
velocity as a function of the number of quadrature points for all
three test cases. 

\begin{figure}
  \centering
  \includegraphics{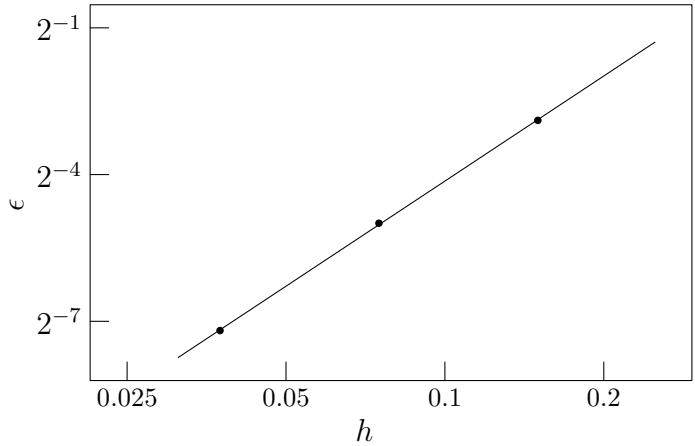}
  \caption{Error in velocity against discretization length $h$, for
    $64\times64$ point quadrature: circles: errors; solid line: linear
    fit $\epsilon=8.03h^{2.15}$}
  \label{fig:convergence:h}
\end{figure}

It is clear that, as proposed in \S\ref{sec:speed}, the accuracy of
the integration is controlled by the order of Gaussian quadrature, up
to a limit fixed by the resolution of the vorticity distribution. This
is investigated in figure~\ref{fig:convergence:h} which shows the
r.m.s error $\epsilon$ of figure~\ref{fig:convergence} plotted against
a typical discretization length scale $h$. The velocities were
computed using~64 quadrature points in $\phi$ and $\theta$. At this
stage the two lower resolution test cases have reached their minimum
error. The error scales as $h^{2.15}$ showing that the method is
approximately second order.

\subsection{Parallelization}
\label{sec:parallel}

\begin{figure}
  \centering
  \includegraphics{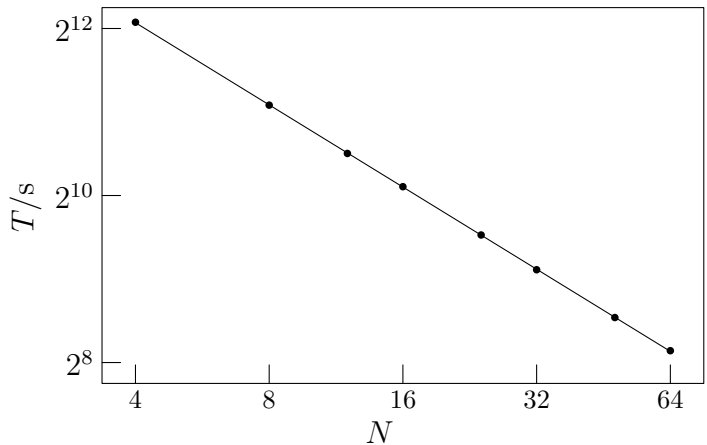}
  \caption{Computation time $T$ for~16000 point Hill's spherical
    vortex as a function of number of processors $N$: circles
    computation time, solid line linear fit.}
  \label{fig:speedup}
\end{figure}

The final requirement of a velocity calculation algorithm is that it
be possible to implement it efficiently on a parallel system. The
method of \S\ref{sec:velocity} has been coded for an MPI system by
sharing the tetrahedra out amongst the processors and combining the
velocity contributions at the end of the calculation. Only one
communication is needed, to combine the velocities computed on each
processor so that the network overhead is low.
Figure~\ref{fig:speedup} shows the computation time for velocity on
Hill's spherical vortex discretized with~16000 points,
with~4--64~3.2\giga\hertz\ processors. The fitted line shows the
calculation time scaling as $N^{-0.98}$, a near linear speedup.

\section{Conclusions}
\label{sec:conclusions}

A method for the evaluation of Biot--Savart integrals has been
presented which eliminates the need for complicated mathematical
operations over most of the calculation. With the exception of a
negligible number of trigonometric factors which must be pre-computed,
the method requires no operation more costly than a division.
Calculations performed using the technique demonstrate that it is
approximately second order accurate, that it converges correctly and
that it achieves near-linear speedup on a parallel system. Future work
will consider how to use the approach in methods similar to fast
multipole and how to improve the ray tracing algorithm used, perhaps
by taking advantage of adjacency relationships between tetrahedra.

In conclusion, we note that the method is based on a ray-tracing
technique from computer graphics and that such methods are now being
implemented at hardware level on mass market systems. This opens the
possibility of using the calculation method on relatively cheap
hardware, giving further acceleration at little additional cost.

\bibliography{abbrev,misc,vortex,turbulence,jets,maths,%
  propnoise,identification}

\end{document}